\documentclass[leqno,12pt]{amsart}
\usepackage[T1]{fontenc}
\usepackage[english]{babel}
\usepackage{amssymb}
\usepackage{amsmath}
\usepackage{amsthm}
\usepackage{amscdm}

\input diagxy.sty
\providecommand{\LyX}{L\kern-.1667em\lower.25em\hbox{Y}\kern-.125emX\@}

\input headers-note.sty
\def\SC{{\mathfrak{SC}}}

\begin{document}
\title[Uniquely divisible abelian extensions of $E(\barq)$ by
 $\Lambda$]{A remark on transitivity of Galois action on the set of
uniquely divisible abelian extensions in $\Ext^1(E(\barq),\Lambda)$}
\author[Misha Gavrilovich]{ Misha Gavrilovich\\
{Balliol College, Oxford}}\address{URL:
http:\!/\!\!/\!misha.uploads.net.ru}
%\maketitle
%\begin{document}
 \abstract{We study Galois action on $\Ext^1(E(\barq),\Z^2)$ and
interpret our results as partially showing that the notion of a path
on a complex elliptic curve $E$ can be characterised algebraically.
The proofs show that our results are just concise reformulations of
Kummer theory for $E$  as well as the description of Galois action
on the Tate module.
%We hope these reformulations provide a new and
%useful point of view on these
%arithmetic results. %Kummer theory of an elliptic curve.
% on as well as the description of Galois action on Tate module.

 Namely, we prove (a),(b) below by showing they are equivalent to (c)
 which is well-known: $(a)$
group $\galqq$ acts transitively on the set of uniquely divisible
abelian $\EndE$-module extensions of $E(\barq)$ of algebraic points
of an
elliptic curve, by $\Lambda\cong\Z^2$, %and interpret this as that
$(b)$~natural algebraic properties characterise uniquely the
Poincar\'e's fundamental groupoid of a complex elliptic curve,
restricted to the algebraic points,
(c) % Kummer theory of the elliptic curve $E$:
 up to finite index, the image of the Galois
action on the sequences $(P_i)_{i>0},jP_{ij}=P_i,i,j>0$ of points
$P_i\in E^k(\barq)$ is as large as possible with respect to linear
relations between the coordinates of the points $P_i$'s.

%We do so by showing that $(a),(b)$ above are essentially
%reformulations of Kummer theory of $E$, as well as the description
%of the image of the Galois representation on the Tate module.

Our original motivations come from model theory.} } \maketitle

{\ \ \ \ \ \ \ \ \ \ \ \ \ \ \ \ \ \ \ \ \ \ \ \ \ \ \ \ \ \ \ \ \ \ \ \ \ \ \ \ \ \ \ \ \ \ \ \ {\tiny%  ...on s'empressera de publier ses moindres
%observations pour peu qu'elles soient nouvelles,
...et on ajoutera:<<je ne sais~pas le reste>> }}

%\epigraph{...on s'empressera de publier ses moindres observations
%pour peu qu'elles soient nouvelles, et on ajoutera: <<je ne sais~pas
%le reste>>}

%\epigraph{Un homme qui a une idée peut choisir entre, avoir, sa vie
%durant, une réputation colossale d'homme savant, ou bien se faire
%une école, se taire, et laisser un grand nom dans l'avenir. Le
%premier cas a lieu s'il pratique son idée sans l'émettre, le second,
%s'il la publie. Il y a un troisième moyen juste milieu entre les
%deux autres, c'est de publier et de pratiquer, alors on est
%ridicule.}

\section{\label{ar:intro}Introduction}

%Motivated by model theory, Zilber~\cite{ZilberCovers} proves that
%$\Aut(\C/\Q)$ acts transitively on the set of uniquely divisible
%abelian group extensions of $\C^*$ by $\Z$, and that in general,
%this holds for any algebraically closed field $K=\bar K$ of zero
%characteristic. He uses model-theoretic tools to treat the case of
%an uncountable field.

\subsection{Abelian group extensions} The universal covering space
of  an elliptic curve $A=E$ is just $\C$; the linear structure on
$\C$ plays an important r\^ole in the theory of (complex) elliptic
curves as algebraic varieties. This suggests that the relevant
linear structure on $\C$ is determined up to isomorphism by the
algebraic curve itself. A uniquely divisible abelian extension of
the group $E(\C)$ with kernel $\Z^2$ describes such a linear
structure; thus, the above considerations suggest that such an
extension is unique, up to an automorphism $\Aut(\C/\Q)$.

%An analysis of
%$\La$ based on these observations shows that $\La$ is a language
%describing the linear structure on $\C$, and the pull-back of the
%algebraic structure on the elliptic curve.

The next proposition partially confirms these suggestions.
%allows us to prove that for $A=E$ an elliptic curve,
%Goal~\ref{asp:atomic} is implied by the following proposition.
We conjecture it holds for $\C$ and in general, any algebraically
closed field of zero characteristic.

\begin{prp-result}\label{prp:ext:0}\label{prp:ext}
Let $E$ be an elliptic curve defined over a number field
$k\subset\barq$. Assume that all the endomorphisms of $E$ are
definable over $k$. If $\EndE=\Z$, then the Galois group $\galqq$
acts transitively on the set of uniquely divisible abelian
extensions $\Ext^1_{\text{AbGroups}}(E(\barq),\Lambda)$ of
$E(\barq)$ by $\Lambda\cong \Z^2$. If $\EndE\neq\Z$, then there are
at most finitely many orbits on the set of  uniquely divisible
$\EndE$-module extensions in
$\Ext^1_{\EndE\text{-}\textrm{mod}}(E(\barq),\Lambda)$ of
$\EndE$-modules.
\end{prp-result}

Here $\Lambda$ denotes the kernel of $\EndE$-module covering map
$\C\ra E(\C)$. We prove the proposition by an inductive argument
using Kummer theory of $E$ and the description of the Galois action
on Tate module; moreover, the proof shows that Proposition~1 is
equivalent to these arithmetic results.

\subsection{An algebraic notion of a path up to homotopy} We ask
whether \emph{the notion of paths (up to fixed point homotopy) on an
elliptic curve $E(\C)$ may be described by its natural algebraic
properties}. ``Paths up to fixed point homotopy'' are usually
thought of in the context of the \emph{Poincar\'e's fundamental
groupoid}, which can be thought of as a 2-functor. Hence, we may
reformulate the question as: \emph{Is the fundamental groupoid
functor on the complex algebraic varieties determined by its natural
algebraic properties up to natural equivalence and an automorphism
of the source category?} Proposition~\ref{prp:functor} is a partial
positive answer to this question.

%Let
%$$\E\subset\factor{\textrm{Var}}{\,\barq},\ \Ob\, \E=\{E^n:n\geq 0\},\
%\Mor_{\!\E}(X,Y)=\Mor_{\text{Var}/\barq}(X,Y)$$ be the full
%subcategory of the category of varieties whose objects are the
%Cartesian powers of $E$ including the trivial Cartesian power
%$E^0=0$. The Galois group $\galqk$ acts on the category
%$\text{Var}/\barq$, and its action restricts to the category $\E$:
%the action leaves each object of $\E$ invariant and permutes the
%morphisms. Note that $\galqk$ acts on the set of functors from $\E$
%via its action on the category $\E$ by automorphisms.

\begin{prp-result}[Universality of fundamental groupoid functor]\label{prp:functor:0}\label{prp:functor}
Let $E$ and a number field $k$ be as above. Let $\E$ be the
full subcategory of $\factor{\textrm{Var}}{\,\barq}$ consisting of
Cartesian powers of $E$.

Let $\uOmega:\E\ra\Groupoids$ be a functor satisfying conditions
(1)-(3) of Definition~\ref{defn:a-fundamental} below (unique
path-lifting along \'etale morphisms, preservation of direct
product, etc).

 Assume further that \bi
\item[(4)] $\uOmega(E)$ is a connected groupoid
\item[(5)] there is an isomorphism
$$\Omega_{0,0}(E):=\{\gamma\in \Omega(E): s(\gamma)=t(\gamma)=0\}\cong \Z^2
\text{ as }\EndE\text{-modules}.$$\ei
 Then, for $E$ with no complex multiplication, there exists an
automorphism $\sigma\in\galqk$ such that the functors
$\uOmega,{\pitop}_{|\barq}\circ\sigma:\E\ra \Groupoids$ are
naturally equivalent:
$$\uOmega\cong {\pitop}_{|\barq}\circ\sigma$$
%\becd \E @>>\uOmega> \Groupoids \\
%   @V\exists\! V\sigma\in\galqk V   @Y\exists\!Y \mathrm{natural\ transformation}Y\\
%            \E @>\pitopq>> \Groupoids \eecd

If $E$ has complex multiplication, then there exists finitely many
functors $\F_1,...,\F_n:\E\ra \Groupoids$ satisfying conditions
(1)-(5) above such that for any functor $\uOmega$ satisfying (1)-(5)
above there exists a Galois automorphism $\sigma\in\galqk$ and a
number $0<i<n+1$ such that $\uOmega$ and $\F_i\circ\sigma$ are
naturally equivalent.
\end{prp-result}

Here $\Pt_{{\pitop}_{|\barq}}(V)=V(\barq)$ and
$$\Omega_{{\pitop}_{|\barq}}(V)=\{\gamma\in{\pitop}(V(\C)):
s(\gamma),t(\gamma)\in V(\barq)\}$$ is the restriction of ${\pitop}$
to $\barq$-rational points.

We prove Proposition~2 by reducing it to Proposition~1.

\subsection{Kummer theory and Galois representations}
We prove Propositions~1~and~2 by a rather simple induction based on
Kummer theory and the image of Galois action on Tate module. The
full proof is carried out in in \S\ref{proof:prp:ext}; here we just
present an argument showing how these arithmetic results may be
relevant. For simplicity assume $\EndE=\Z$ and consider an extension
$H\ra^\phi E(\barq)$.  Take a pair $\lambda_1,\lambda_2\in \Lambda$
of generators of abelian group $\Lambda=\ker\phi$;  the extension
then gives rise to the infinite sequences
$(\phi(\lambda_1/n))_n,(\phi(\lambda_2/n))_n$ which can be thought
of as elements of Tate module $\prod_l T_l(E)$. If we take another
extension $H'\ra^{\phi'} E(\barq)$, we can obtain similar sequences
$((\phi(\lambda'_1/n))_n,(\phi(\lambda'_2/n))_n$. If the extensions
$H\ra^\phi E(\barq)$ and $H'\ra^{\phi'\circ\sigma} E(\barq)$ are
isomorphic for some $\sigma\in\galqk$ and an isomorphism $':H\ra H'$
takes $(\lambda_1,\lambda_2)$ to $(\lambda'_1,\lambda'_2)$, then
these pairs of sequences are conjugated by Galois action. Were
$\galqk=\SL_2(\hz)$, it would be so
%these sequences would be conjugated
for any two linearly independent pairs $(\lambda_1,\lambda_2)$ and
$(\lambda'_1,\lambda'_2)$ of generators.

If  $\galqk$ is of finite index in $\SL_2(\hz)$, then it is
sufficient to take $\lambda'_1,\lambda'_2$ such that
$\phi'(\lambda'_1/N)=\phi(\lambda_1/N),
\phi'(\lambda'_2/N)=\phi(\lambda_2/N)$, for some $N$ large enough;
that is possible due to surjectivity of $\phi'$ and the unique
divisibility of $H'$.

\subsection{Model theory: motivations and generalisations} The
results of this note appear naturally in a model-theoretic framework
of <<logically perfect structures>> developed by Zilber
\cite{ZilberLP}; in this way they have been obtained in the author's
DPhil thesis \cite{gavr-thesis}. We believe that the context and
techniques of model theory are essential to generalise our results
to other varieties and fields of arbitrary cardinality including the
field of complex numbers. However, in this note we do not discuss
the model theoretic motivations and techniques; we refer to the the
author's DPhil thesis \cite{gavr-thesis} for that. The note presents
only those results of the thesis which can be stated and proven in a
model-theory free fashion; in particular, {\em the note omits the
discussion of the Shafarevich conjecture on holomorphic convexity of
universal covering spaces and the case semi-Abelian varieties of
higher dimension.} The results follow a line originated by Zilber
\cite{ZilberCovers} treating fully the case of the multiplicative
group of an algebraically closed field of characteristic 0.

 We would appreciate references to any literature discussing similar
questions.

\subsection{Structure of the paper}

%In Proposition~\ref{prp:ext:0} of \S\ref{1.1} we state the main
%observation, also including the case of an elliptic curve with
%complex multiplication. In
% Proposition~\ref{prp:functor:0} of \S\ref{1.2} we state the universality property of
% the Poincar\'e fundamental groupoid functor.

%restricted to the full subcategory of Cartesian powers of an
%elliptic curve. What we would like to do is define a natural class
%of functors such that for any functor $\F$ of the class, there
%exists $\sigma\in\Aut(\C/\Q)$ such that $\pi_1\cong \F\circ\sigma$
%are naturally equivalent (we think of $\sigma$ as an automorphism of
%the category of varieties over $\C$); in fact we would like to prove
%this for $K=\bar K$ an algebraically closed field of any
%cardinality. Proposition~\ref{prp:functor} claims the universality
%property for $K=\barq$ the field of algebraic numbers.

We state our results in detail in \S\ref{iv:results}.

We prove Proposition~\ref{prp:ext:0} in \S\ref{proof:prp:ext}. We
establish the equivalence of the algebraic approach of
Proposition~\ref{prp:ext} and the topological approach of
Proposition~\ref{prp:functor:0} in \S\ref{sec:functor}. We state a
precise conjecture about Shimura curves in \S\ref{sec:shimura}.

\section{\label{iv:results}Results}

In this section we state our results in full, and hint on a
connection between the reformulations.

\subsection{\label{1.1}Uniquely divisible extensions of Abelian groups}

Let $E$ be an elliptic curve defined over a number field
$k\subset\barq\subset\C$, and let \becd 0 @>>> \Lambda @>>> \C @>p>>
E(\C) @>>> 0. \eecd be the universal  covering of  $E(\C)$. Let
$0\in E(k)$ denote a $k$-rational point which is zero of the
additive group $E(\C)$.

\begin{prp}\label{prp:ext:0}\label{prp:ext}  Assume
that all the endomorphisms of $E$ are definable over $k$. If
$\EndE=\Z$, then the Galois group $\galqq$ acts transitively on the
set of uniquely divisible abelian extensions
$\Ext^1_{\text{AbGroups}}(E(\barq),\Lambda)$ of $E(\barq)$ by
$\Lambda\cong \Z^2$. If $\EndE\neq\Z$, then there are at most
finitely many orbits in uniquely divisible extensions in
$\Ext^1_{\EndE\text{-}\textrm{mod}}(E(\barq),\Lambda)$ of
$\EndE$-modules.
\end{prp}
%\bi
%\item There exists a uniquely divisible extension
%in $\Ext^1_{\EndE\text{-mod}}(E(\barq),\Lambda)$.
%\item\bi
%\item if $\EndE=\Z$, then the Galois group acts transitively on the
%set of uniquely divisible extensions
%$\Ext^1_{\text{AbGroups}}(E(\barq),\Lambda)$.
%\item the action of the Galois group has only finitely many orbits
%on the set of  uniquely divisible extensions in
%$\Ext^1_{\EndE\text{-mod}}(E(\barq),\Lambda)$. \ei \ei
%\end{prp}
%
%An isomorphism $h:V\ra W$, $h(\Lambda)=\Lambda$ induces an
%$\EndE$-linear map $\tau:E(\barq)\ra E(\barq)$; the
%proposition~(item~(a)) says we may modify $h$ to have $\tau$ induced
%by a Galois automorphism. In fact, we show we can extend any
%$h,\sigma$ defined on \emph{finite} dimensional submodule of $V$ to
%a total map $h:V\ra W,\sigma\in\galqq$.
% The proposition
%may seem somewhat surprising: say for elliptic curves without
%complex multiplication, for every automorphism $\tau:E(\barq)\ra
%E(\barq)$ as \emph{an Abelian group}, for the $\phi:V \ra E(\barq)$
%and $\tau\circ\phi:V \ra E(\barq)$ there exists a Galois
%automorphism $\sigma\in \galqk$ and an automorphism $\xi\in \Aut(V)$
%such that $\phi\circ\xi=\sigma\circ\tau\circ \phi$.

We conjecture the proposition holds for any algebraically closed
field of zero characteristic, i.e.~the action of $\Aut(\bar K/\Q)$
on the set of uniquely divisible abelian extensions in
$\Ext^1(E(\bar K),\Lambda)$. Zilber~\cite{ZilberCovers} proves the
transitivity of the $\Aut(\bar K/\Q)$-action on $\Ext^1(\bar
K^*,\Z)$ for arbitrary algebraically closed field $K=\bar K$ of
characteristic 0.

Remark that the set of non-equivalent extensions
$\Ext^1_{\EndE\text{-mod}}(E(\barq),\Lambda)$ is of cardinality
$2^{\aleph_0}$. Also note that the injectivity of the profinite
completion $\widehat\Lambda$ of the kernel $\Lambda$ implies
$\Ext^1_{\widehat\EndE\text{-mod}}(E(\barq),\widehat\Lambda)=0$
making the proposition trivially true in this case.

A way to think of it is that the proposition claims that \emph{it is
possible to describe the universal covering space of an elliptic
curve in a purely algebraic way}, admittedly with respect to a
rather weak, linear structure on it.

%\begin{rem} We believe that our treatment of elliptic curves with
%complex multiplication is incomplete; the formulation of $2(b)$ that
%an extension is equivalent to one of finitely many is ugly. It
%should be possible to prove a full analogue of $1(b)$, after perhaps
%some suitable reformulation. This is particularly apparent as there
%is a complete description of Galois action on $E_\mathrm{tors}$
%provided by the theory of complex multiplication on elliptic
%curves.\end{rem}

\subsection{\label{1.2}Relation between Propositions~1~and~2}

Here we interpret Proposition~\ref{prp:ext} as a universality
property of the Poincar\'e fundamental groupoid functor. Identify
points in $\C$ and homotopy classes of paths in $E(\C)$ starting at
$1$ via the period map $\gamma\longmapsto\int\limits_\gamma dz.$
Then the addition on $\C$ corresponds to point-wise addition of
paths, and dividing by $n$ corresponds to path-lifting along the
\'etale morphism $nx:E(\C)\ra E(\C)$. Using this correspondence, we
interpret the unique divisibility of $\C$ as the unique path-lifting
property along \'etale morphisms from $E(\C)$ to $E(\C)$.
%The formula
%$$\int\limits_{\gamma_1} dz+\int\limits_{\gamma_2} dz=\int\limits_{\gamma_1\cdot\gamma_2} dz=
%\int\limits_{\gamma_1\gamma_2} dz$$ shows that the addition in $\C$
%corresponds to the point-wise multiplication of paths. Similarly,
%unique divisibility by a natural number $n$ corresponds to the
%unique path-lifting along \'etale morphism $x^n:\Cm\to\Cm$.

These observations allow us to reformulate the previous result as an
algebraic characterisation of the Poincar\'e's fundamental groupoid
of a complex elliptic curve, i.e.~that \emph{the notion of paths (up
to fixed point homotopy) on a complex elliptic curve $E(\C)$ may be
described
by its natural \texttt{algebraic} properties}. %``Paths up to fixed point
%homotopy'' are usually thought of in the context of the
%\emph{Poincar\'e's fundamental groupoid}. %which can be thought of as
%a 2-functor.

A precise reformulation of the question is: \emph{Is the fundamental
groupoid functor on the complex algebraic varieties determined by
its natural algebraic properties up to natural equivalence and an
automorphism of the source category?}
Proposition~\ref{prp:functor:0} provides a partial positive answer
to this question; the proposition is just an equivalent
reformulation of Proposition~\ref{prp:ext}.

\subsection{\label{sec:two-functor}The universality property of the Poincar\'e's
fundamental groupoid functor}

In this subsection we introduce in detail the notions of category
theory appropriate to state Proposition~2.

\subsubsection{The example: fundamental groupoid functor ${\pitop}(E(\C))$ as a two-functor}

It is  convenient to consider 2-functors instead of functors to
groupoids; the notions are equivalent. Before defining a 2-functor
formally, we illustrate the notion by an
example of Grothendieck~\cite{grothendieck-quillen}. % which have motivated the notion.

The path 2-functor  $\uOmega$ on the category $\Top$ of topological
spaces is a tuple $(\Pt,\Omega,s,t,\cdot)$ consisting of a
\emph{functor of points} $\Pt:\Top\ra\Sets$ and a \emph{paths
functor} $\Omega:\Top\ra\Sets$ together with the following data:
       \bi
       \item $\Pt(T)$ is the set of points of topological
       space $T$ and the morphism $\Pt(f):\Pt(T_1)\ra \Pt(T_2)$ is $f:T_1\ra T_2$ as a map of sets.
       \item \label{omega} $\Omega(T)$ is the set of all paths in topological
       space $T$, i.e.~continuous functions $\gamma:[a,b]\ra T$, $a,b\in \R$; similarly
       $\Omega(f),f\in\Hom_{\Top}(T_1,T_2)$ is the map  taking a path $\gamma:[a,b]\ra T_1$ into
       $f\circ\gamma:[a,b]\ra T_2$.
       \item $s_T,t_T:\Omega(T)\ra \Pt(T)$ are functions from the set of
       paths in $T$ to their endpoints in $T$; the function
       $s(\gamma)=\gamma(a)$ (\emph{source}) gives the beginning point of a path $\gamma$,
       and the function $t(\gamma)=\gamma(b)$ ({\em target}) gives the ending
       point of path $\gamma$.
       \item $\cdot_T:\Omega(T)\times\Omega(T)\ra\Omega(T)$ is the partial operation of
       the concatenations of paths, taking $\gamma_1:[a,b]\ra T$,
       $\gamma_2:[b,c]\ra T$ into $\gamma=\gamma_1\gamma_2$,
       $\gamma_{|[a,b]}=\gamma_1,\gamma_{|[b,c]}=\gamma_2$.
\ei

Thus, a 2-functor from $\Top$ to $\Sets$ consists of two functors
$\Pt,\Omega:\Top \ra \Sets$, and two natural transformations
$s,t:\Omega \ra \Pt $ from functor $\Omega$ to $\Pt$; $s$ stands for
source and $t$ stands for target. For each $T$, there is also a
functorial associative operation $\cdot_T$ defined on
$\Omega_{x,y}(T)\times \Omega_{y,z}(T)\ra \Omega_{x,z}(T)$, where
$\Omega_{x,y}(T)=\{\gamma\in\Omega(T): s(\gamma)=x, t(\gamma)=y\}$,
etc; the operation $\cdot$ makes $\Omega_{x,x}(T)$ into a group; in
the example above, $\Omega_{x,x}(T)$ is the set of all loops in $T$
based at $x\in T$.

In particular, for each $T$ the set $\Omega(T)$ carries the
structure of a groupoid; in fact, it is conventional to consider
$\uOmega$ as a functor to the category of groupoids.

If in item $(\ref{omega})$ we define $\Omega(T)$ to be the set of
all paths \emph{up to fixed point homotopy}, then we obtain the
notion of the fundamental groupoid functor. The advantage of the
original definition is that one may try and define $n$-functors
describing $n$-dimensional homotopies on topological space $T$,
cf.~Grothendieck~\cite{grothendieck-quillen} for motivations;
Voevodsky-Kapranov~\cite{voevodsky-kapranov}  propose %????????
an exact definition. Grothendieck~\cite{grothendieck-quillen}
explains that it is essential not to insist on strict associativity
etc, but rather to consider {\em all the identities to hold up to a
homotopy of higher dimension}. This may be useful to generalise
Proposition~2.

\subsubsection{Abstract fundamental groupoid functors}

Here we define the path-lifting property for a 2-functor, and an
abstract fundamental groupoid functor as a functor preserving direct
products, possessing the path-lifting property and with a particular
functor of points.

\begin{defn}\label{defn:a-fundamental} Let $\uOmega$ be a 2-functor from a subcategory $\E$ of
the category of varieties over an algebraically closed field $K$,
into $\Sets$. We say $\uOmega$ is an \emph{abstract fundamental
groupoid $K$-valued functor} if $\uOmega$ satisfies the following
properties:
\begin{enumerate}\label{functor:cond}
\item the functor $\fun$ is the functor of $K$-rational points:
\[ \fun(X)=X(K)=\Mor_{\mathrm{Var}/K}(0,X)\]
(here $0$ denotes a single point variety defined over $k$)
\item the functor $\path$ preserves the
direct product:\becd \path(X\times Y)=\path(X)\times\path(Y)\\
                \path(f\times g)= \path(f)\times \path(g)\\
                s_{X\times Y}=s_X\times s_Y, t_{X\times Y}=t_X\times t_Y\\
                (\gamma_1\times
                \gamma_2)\cdot_X(\gamma'_1\times\gamma'_2)=
                (\gamma_1\cdot_X\gamma'_1)\times(\gamma_2\cdot_X\gamma'_2)
                \eecd (the direct product is taken in the category of $\Sets$.)
\item {unique path-lifting property:} if $p\in\Hom_{\E}(X,Y)$ is an \etale{} morphism of algebraic
varieties, then for any points $x,y\in \Pt(X)$ the map
\[
\path(p):\bigcup_{y\in\fun(X)}\Omega_{x,y}(X)\ra\bigcup_{z\in\fun(Y)}\Omega_{p(x),z}(Y)\]
 is a bijection.
%\item {van Kampen theorem}: for any $X$ and $\gamma\in\Fun(X),$
%if $X=X_{1}\cup\ldots\cup X_{n}$ then there exists a number $N$
%and a sequence $\gamma_{i},0<i<N$ such that \[
%\gamma=\gamma_{1}\cdot\ldots\cdot\gamma_{N}\] and \emph{\[
%\gamma_{i}\in\Fun(X_{1})\cup\ldots\cup\Fun(X_{N})\] }
\end{enumerate}
\end{defn}

\subsubsection{Universality of the fundamental groupoid functor}

Let
$$\E\subset\V,\ \Ob\,\E=\{E^n:n\geq0\},\ \Mor_\E(X,Y)=\Mor_{\mathrm{Var}/K}(X,Y)$$
be the full subcategory of the category of varieties whose objects
are the Cartesian powers of $E$ including $E^0=0$ a variety
consisting of the single point $0$. By the definition of a full
subcategory, the morphisms of $\E$ are morphisms of varieties
between the objects of $\E$.

The Galois group $\galqk$ acts on the category $\E$; the action
leaves the objects invariant but permutes the morphisms. % not necessarily preserving $0\in E(k)$;
Recall we assume that all endomorphisms of $E$ preserving $0\in
E(k)$ are defined over its field $k$ of definition.

%Note that $\galqk$ acts on the set of 2-functors from $\E$: it acts
%on the functor of points via the formula
%$\Pt(V)=\sigma(\Pt(V))=\sigma(V(\barq))$ while leaving the functor
%of paths $\Omega $ and functions $s_V,t_V$ invariant.

Recall a groupoid $\Omega(E)$ is \emph{connected} iff  for every
$x,y\in\Pt(E)$ there exists $\gamma\in \Omega(E)$ ``going from point
$x$ to point $y$'', i.e.~$x=s(\gamma),y=t(\gamma)$. The set of such
paths satisfying $x=s(\gamma),y=t(\gamma)$ is denoted by
$\Omega_{x,y}(E)$.

Using the notion of an abstract fundamental groupoid functor, we
restate Proposition~\ref{prp:functor}; then we show it is
essentially equivalent to Proposition~\ref{prp:ext}. The proof
basically reconstructs the ``universal covering space'' $V$ as the
set of all paths $\bigcup_{*\in E(\barq)}\Omega_{0,*}(E)$ leaving a
particular point; functoriality of $\uOmega$ allows us to define
$\EndE$-module structure on $V$; the unique path-lifting property of
$\Omega$ ensures unique divisibility.

\begin{prp}[Universality of fundamental group functor]\label{prp:functor:2}
Assume that $E$ is  an elliptic curve defined over a number field
$k$ and does not have complex multiplication. Let $\E$ be the full
subcategory of Cartesian powers of $E$ as above.

Let $\uOmega$ be an abstract fundamental groupoid $\barq$-valued
functor on category $\E$. Assume  \bi \item $\uOmega(E)$ is a
connected groupoid
\item there is an isomorphism
$$\Omega_{0,0}(E)\cong \Z^2
\text{ as }\EndE\text{-modules}.$$ \ei Then there exists an
automorphism $\sigma\in\galqk$ such that the 2-functors
$\uOmega,{\pitop}_{|\barq}\circ\sigma:\E\ra \Sets$ are naturally
equivalent (as 2-functors),
$$\uOmega\cong{\pitop}_{|\barq}\circ\sigma$$
%pçppl
  \becd \E @>>\uOmega> \Sets \\
       @V\sigma VV   @YY\mathrm{natural\ transformation}Y\\
  \E @>\pitopq>> \Sets \eecd

If $E$ has complex multiplication, then there exists finitely many
abstract fundamental groupoid $\barq$-valued functors
$\F_1,...,\F_n:\E\ra \Groupoids$ such that for any other abstract
fundamental groupoid $\barq$-valued functor $\uOmega$ there exists a
Galois automorphism $\sigma\in\galqk$ and a number $0<i<n+1$ such
that $\uOmega$ and $\F_i\circ\sigma$ are naturally equivalent.

\end{prp}

Recall $\Pt_{{\pitop}_{|\barq}}(V)=V(\barq)$ and
$\Omega_{{\pitop}_{|\barq}}(V)=\{\gamma\in{\pitop}(V(\C)):
s(\gamma),t(\gamma)\in V(\barq)\}$.

%The reduction of holds for the general case of an elliptic curve.

The condition in the definition of an abstract fundamental groupoid
functor are somewhat reminiscent of the conditions defining the
scheme-theoretic algebraic fundamental group $\pi_1^\textrm{alg}$
[SGA$4\frac 12$]; however, there $\pi_1^\textrm{alg}$ takes values
in the category of \emph{profinite} groups, in particular
$\pi_1^\textrm{alg}(\Cm,1)=\hz$,
$\pi_1^\textrm{alg}(E(\C),0)=\hz^2$.
%In terms of model theory, this
%corresponds to considering \emph{homogeneous models}.

It is  natural to consider whether a path lies in an algebraic
subvariety, so if $\uOmega$ is provide a useful notion of a path on
$E(\C)$,
%thus our notion of a path given via
the 2-functor $\uOmega$ \emph{restricted} to $\E$ should be able to
express when (a representative of the homotopy class of) a path lies
in an {\em arbitrary} algebraic subvariety. This is indeed the case:

\begin{rem}[Recovering $\Omega(Z)$ for arbitrary closed subvariety $Z$ of
$E^n$]\label{rem:OmegaPaths} The information contained in the
functor ${\pitop}|\E$ restricted to the full subcategory $\E$ of
Cartesian powers of an elliptic curve is enough to determine whether
a path lies in a closed subvariety. The key fact here is that for a
normal subgroup $H\vartriangleleft\pi_1(E^n(\C))$, there exists an
\emph{H-Shafarevich} morphism $\mathrm{Sh}_H:E^n\ra E^m$ such that
for an arbitrary irreducible $Z\subset E^n(\C)$, it holds $Z\subset
\ker f$ iff the image $\im[\pi_1(\hat Z,z)\ra \pi_1(E^n(\C),z)]$ has
a finite index subgroup contained in $H$. % Lemma~\ref{lem:simcBZ}
%states this in a rather different language; cf.~also
%Remark~\ref{rem:alb} for a technical explanation of the connection
%between the reformulations.
\end{rem}

%A homotopy class of paths has a representative lying  in an open
%algebraic subvariety iff it has one lying in its closure; this holds
%at least for normal varieties. An extension of this observation
%based on Fact~\ref{fact:piiopen} and considerations
%in~\S\ref{sec:nonnormal} shows that the functor $\uOmega$ is enough
%to determine whether the homotopy class of a path has a
%representative lying in an arbitrary constructible subset.

\section{\label{proof:prp:ext}Extensions of $E(\barq)$ by $\Lambda$}

In this section we state and prove Proposition~\ref{prp:functor};
the proof is a model theoretic argument based on Kummer theory and
the description of the image of Galois representation on Tate module
$T(E)$. However, we tried to be very explicit and have avoided any
model-theoretic terminology in the exposition of the proof. The only
model theory left in the proof is in the level of motivations; as
was noted earlier, those model-theoretic motivations are useful to
try and generalise the formulation to other varieties and contexts
in general.

\subsection{Transitivity of $\galkq$ action on the uniquely divisible
$\End E$-module extensions of $E(\barq)$ by $\Lambda$}

To fix notations for the proof, we restate Proposition~\ref{prp:ext}
in an expanded form.

\begin{prp}\label{prp:ext:1} Let $E(\C),\Lambda$ be as above.

 \bi\item \label{p:exist} There exists a uniquely divisible $\EndE$-module $V$
 and a short exact sequence of $\EndE$-modules
  \becd 0 @>>> \Lambda @>>> V @>>> E(\barq) @>>> 0. \eecd

 \item \label{p:iso} \bi
 \item If $E$ has no complex multiplication and $\EndE=\Z$, then for any
 uniquely divisible $\Z$-module extensions $W,V$
 of $E(\barq)$ by $\Lambda$,
 fitting into the short exact sequences as above, there exist a commutative diagram:
  \becd\label{fml:exact:ell:a} 0 @>>> \Lambda @>>> V @>\phi>> E(\barq) @>>> 0 \\
  @XXX @U {\exists h_{|\Lambda}} UU @U\exists\! U {\!h\in\Hom(V,W)} U  @UU \exists\sigma\in\galqk U   \\
  0 @>>> \Lambda @>>> W @>\psi>> E(\barq) @>>> 0 \eecd
  \item
  If $\EndE\neq\Z$, then there exist finitely many uniquely divisible $\EndE$-module
  extensions $W_1,...,W_n$  of $E(\barq)$ by $\Lambda$,
 fitting into the short exact sequences as above, such that for any uniquely
 divisible $\EndE$-module extension $V$ there exist a commutative diagram:
  \becd\label{fml:exact:ell:2} 0 @>>> \Lambda @>>> V @>\phi>> E(\barq) @>>> 0 \\
  @XXX @U {\exists h_{|\Lambda}} UU @U\exists\! U {\!h\in\Hom(V,W)} U  @UU \exists\sigma\in\galqk U   \\
  0 @>>> \Lambda @>>> W_i @>\psi_i>> E(\barq) @>>> 0 \eecd\ei\ei
\end{prp}

\bp %By \cite{silvermanAEC}, there exists an embedding $k\hra\C$ of
%the field $k$ of definition of $E$ such that the kernel of the
%universal covering map is isomorphic to any such $\Lambda$:
By assumption there is a covering
 \becd 0 @>>> \Lambda @>>> \C @>p>>
E(\C) @>>> 0 \eecd

The endomorphisms act on the complex plane $\C$ by multiplication by
complex numbers, and so in particular $\C$ is a uniquely divisible
$\EndE$-module. The set $E(\barq)$ of points of $E$ over an
algebraically closed subfield is closed under addition and
$\EndE$-multiplication i.e.~is an $\EndE$-submodule, necessarily
uniquely divisible; so is then $V=\pinv(E(\barq))$; take that to
obtain a short exact sequence as above. This proves
$(\ref{p:exist})$.

 To prove $(\ref{p:iso})$, we use the Kummer theory of elliptic
 curves and the Serre's results
about the image of Galois representation on the torsion points.
Essentially, what we need is that all the restrictions on the image
of the Galois action on the Tate module $\prod_lT(E)$ come from the
geometry of $E$ and in this case, are linear. That is exactly what
is provided by the Kummer theory and other results we use.

Pick a maximal linearly independent set $v_0,v_1,v_2,..\in V$; let
$V_n=\EndE v_0+\ldots+\EndE v_n$ be the submodule generated by
$v_0,\ldots,v_n$, and let $\qV_n=(\EndE)\inv V_n=\{v:\exists
N\in\N(Nv\in V_n)\}$ be its divisible closure. We construct by
induction a partial $\EndE$-module linear map $h_n:\qV_n\ra W_i$
inducing a partial Galois map $\sigma_n:\phi(\qV)\ra E(\barq)$ so
that $h=\cup h_n$ is an isomorphism of $V$ and $W_i$; then the
construction implies $\sigma=\cup\sigma_n$ is a total Galois map on
$E(\barq)$ to $E(\barq)$, and thus there is a commutative diagram as
above.

For a point $a\in E(\C)$, let us call \emph{a compatible sequence of
division points starting at $a$} a sequence $(\aaa {} i)_{i\in\N} $
satisfying $\aaa {} i=j\aaa {} {ij},i,j\in\N$. Compatible sequences
of division points starting at $0\in E(\barq)$ form a \emph{Tate
module $T(E)$.} A compatible sequence of \emph{$l$-primary} division
points is a subsequence $(\aaa {} {l^j})_{j\in \N}$ of a compatible
sequence of division points; such sequences of $l$-primary points
starting at $0$ form {\em $l$-adic Tate module} $T_l(E)$.

\subsection{\label{sec:galois-image}The image of Galois representations on Tate module
$T_l(E)$}
\subsubsection{Base of induction}

Choose $v_0\in \Lambda=\ker\phi$; then $\Q V_0= \phi\inv(\Etors)$,
where $\Etors=\{x\in E(\barq): nx=0 \text{ for some }n\in \N\}$ is
the set of \emph{torsion} points of $E(\barq)$.

As the base step of the induction, in this subsection we construct a
commutative diagram:
  \belcd\label{fml:exact:tors}
  0 @>>> \Lambda @>>> V_0 @>\phi>> \Etors @>>> 0 \\
  @XXX @U {h_0} UU @UU {\!\!\!\!\!\!h_0\in\Hom(V_0,W)} U  @UU \sigma\in\galqk U   \\
  0 @>>> \Lambda @>>> W @>\psi\circ\tau_i>> \Etors @>>> 0 \eelcd
where $\tau_1,...,\tau_n$ is some fixed finite collection of
$\EndE$-endomorphisms of $\Etors$ independent of $V,W$. For curves
without complex multiplication, there is no need to consider
$\tau_i$'s, i.e.~$n=1$, $\tau_1=\id$.

Note that when $E$ has complex multiplication, the diagram above is
somewhat reminiscent of a diagram appearing in the main theorem of
complex multiplication \cite[Chapter 5, Theorem 5.4]{Shimura},
also Lang \cite{lang-complex-multiplication}. % and Rubin\cite{rubin-complex-multiplication-notes}.
Global class field theory provides an explicit description on the
image of Galois action in $\Aut(\Etors)$ in terms of the action of
the ring of  id\'eles on some lattices in the universal covering
space.

\subsubsection{\label{sec:cm}$E$ has complex multiplication}

Pick an arbitrary isomorphism $h_0:\Lambda\ra \Lambda$ and extend it
uniquely to $h_0:V_0\ra W$; we may do so by unique divisibility of
$V$ and $W$. Define an $\EndE$-morphism $\tau:\Etors\ra\Etors$ by
$\tau(x)=\psi\circ h_0\circ \phi\inv(x)$. The calculation $\psi\circ
h_0 (y+\Lambda)=\psi\circ h_0 (y)$ shows it is well-defined;
linearity of $\tau$ follows from that of $\phi,h_0$ and $\phi$.

Ideally we would like to be able to choose $h_0:\Lambda\ra\Lambda$
so that $\tau=\tau_h:\Etors\ra\Etors$ is induced by a Galois
automorphism. Here we prove a weaker statement below.

Denote $\Oo=\EndE$ and $E[n]=\{x\in E(\barq): nx=0\}$ the
$n$-torsion of $E$. The set $E[n]$ is a free 1-dimensional $\factor
\Oo{n\Oo}$-module (\cite[Ch.8\S15,Fact 1]{LangDA}), and
$\Aut_\Oo(E[n])=\Aut_{\Oo\!/ {n\Oo}\mathrm{-mod}}(E[n])\cong
(\factor\Oo{n\Oo})^*$, where $(\factor\Oo{n\Oo})^*$ denotes the
group of invertible elements in $(\factor\Oo{n\Oo})^*$.

An $\Oo$-automorphism of $\Etors$ is given by a compatible system of
$\Oo$-automorphisms of $E[n]$, $n>0$; thus we see that there is an
action of  $\hat\Oo=\lim\limits_n\factor\Oo{n\Oo}$ on $\Etors$ as an
$\EndE$-module; the fact that $\Aut_\Oo(E[n])=\Aut_{\Oo\!/
{n\Oo}\mathrm{-mod}}(E[n])\cong (\factor\Oo{n\Oo})^*$ implies that
$\Aut_\Oo (\Etors)\cong \hat\Oo$.

Now we refer to a consequence of the main theorem of complex
multiplication, namely that, in notation of \cite[Ch.8\S15,Fact
2]{LangDA}, the image of Galois group
\[\label{fml:cmimage} G_{K}=Gal(K(\Etors):K)\ra\prod((\EndE)_{l})^{*}\] is open of finite index,
i.e.~$\im G_K $ is a finite index subgroup of
$\hat\Oo^*=\prod_l\Oo_l^*$. Choose $\tau_1,...,\tau_n$ to be
representatives of conjugacy classes $\factor {\Oo^*\!}{\,\im G_K}$;
we then have that for some $i$ $\tau_i\tau=\sigma\in \im G_K$; this
choice of $h_0,\sigma=\tau_i\tau$ makes the diagram
$(\ref{fml:exact:tors})$ commutative, as required.

\subsubsection{\label{sec:nocmimage}$E$ does not have complex multiplication}

Assume that $E$ does not have complex multiplication, i.e.~$\End
E\cong\Z$; identify $\EndE=\Z$, $T(E)=\hz^2$, and
$\Aut(T(E))=\SL_2(\hz)$. The maps $\phi:V\ra E(\barq)$, $\psi:W\ra
E(\barq)$ define  embeddings $\ita_\phi:\Lambda\ra T(E)$,
$\ita_\psi:\Lambda\ra T(E)$ by
$$\ita_\phi: \lambda\mapsto (\phi(\lambda/j))_{j\in\N},$$
$$\ita_\psi: \lambda\mapsto (\psi(\lambda/j))_{j\in\N}.$$

The images $\ita_\phi(\Lambda),\ita_\psi(\Lambda)$ of the both maps
are dense in $T(E)$ due to the surjectivity of
$\phi,\psi:\Q\Lambda\ra \Etors$.

Take a pair of elements $\lambda_0,\lambda_1\in \ker\phi\cong
\Lambda$ generating $\Lambda$ as an Abelian group; we want to find
$\lambda'_0,\lambda'_1\in \ker\psi\cong\Lambda$ and
$\sigma=\sigma_0\in\galqk$ such that
$$\sigma\ita_\phi(\lambda_0)=\ita_\psi(\lambda'_0),$$
$$\sigma\ita_\phi(\lambda_1)=\ita_\psi(\lambda'_1).$$

Under identification $\ker\phi=\Z^2$, since vectors
$\lambda_0,\lambda_1\in \Z^2$ generate lattice $\Z^2$,  it holds
that $\det(\lambda_0,\lambda_1)=1$. That implies that
$\det(\ita_\phi(\lambda_0),\ita_\phi(\lambda_1))$ is a unit in
$\hz=\prod_l \Z_l$. Similarly
$\det(\ita_\psi(\lambda'_0),\ita_\psi(\lambda'_1))$ has to be a unit
in $\hz$.

By \cite[Theorem 3]{bertrand},  the image of the Galois group
$\galqk$ in the automorphism group $\Aut(T(E))=\SL_2(\hz)$ contains
an open subgroup $\SL_2(N\hz)=\ker\left(\SL_2(\hz)\ra \SL_2(\factor
\Z {N\Z})\right)$, for some $N\in\N$ large enough. Now, there is
$L\in\SL_2(N\hz)$ such that $L(u_0)=u'_0$ and $L(u_1)=u'_1$ iff
$\det(u_0,u_1)=\det(u'_0,u'_1)$ and $u_0\in u'_0+ N \hz$, $u_1=u'_1+
N\hz$.  Thus it is enough to take $\Z$-linearly independent
$\lambda'_0,\lambda'_1$ such that
$\psi(\lambda'_0/N)=\phi(\lambda_0/N)$ and
$\psi(\lambda'_1/N)=\phi(\lambda_1/N)$. Necessarily then
$\det(\ita_\psi(\lambda'_0),\ita_\psi(\lambda'_1))$ is a unit. This
implies that $\det(\lambda'_0,\lambda'_1)=1$ and thus,
$\lambda'_0,\lambda'_1$ generate $\Z^2$. The latter statement is
independent of the identification $\ker\psi=\Z^2$, and this
concludes the proof in the case of no complex multiplication.

%\begin{rem} Note that if the image of Galois were smaller then $\SL_2(\Z)$,
%the argument would fail. This is important, as for in higher
%dimensions $\galqk$ has to stabilise the Weil pairing which is a
%form on the Tate module $T(A)$ of the Abelian variety $A$. This
%implies the straightforward generalisation does not work. However,
%see Remark~\ref{rem:chern} for more on Weil pairing and
%generalisations.\end{rem}

\subsection{Kummer theory}

\subsubsection{The main statement of the Kummer theory for an elliptic
curve}

Let us state the main lemma in a form convenient to us to make an
inductive process; that is a form natural  from the model-theoretic
point of view and corresponds to the property of \emph{atomicity of
certain formulae over the kernel}.
%In the model-theory language, it
%says that the type of a linearly independent tuple in $\uU$ is
%isolated by a quantifier-free formula of language $\Lell$.

\begin{lem}[Kummer theory for an elliptic curve]
\label{lem:kummer:ell} Let $E$ be an elliptic curve defined over a
number field $k$.
 Let
$a_1,\ldots,a_n\in E(\barq)$ be a sequence of points linearly
independent over $\EndE$. Then there exists $N\in\Z$ such that any
two compatible sequences $(\aaa 1 i , \ldots, \aaa n i)_{i\in\N}$,
$(\bbb 1 i , \ldots, \bbb n i)_{i\in\N}$  of division points in
$E(\barq)$ starting at $a_1,\ldots,a_n$ and such that $\aaa 1 N=\bbb
1 N, \ldots, \aaa n N =\bbb n N$, are $\galqk$-conjugated by
$\sigma\in \galqk$,
$$\sigma(\aaa i j) = \bbb i j, \text{ for all }0\leq i\leq n,j\in \N.$$
\end{lem}
\bp[Proof of Lemma] See Bashmakov \cite{bashmakov} for original
results for elliptic curves; see
\cite{ribet-kummer-theory-of-Abelian-tori,Ribet,RibetCohom} for
Kummer theory of Abelian varieties, and see \cite{bertrand} for a
summary of results of Kummer theory of Abelian varieties; we quote
\cite[Theorem 2]{bertrand} from that paper.

We now introduce the notations of  Bertrand~\cite[Theorem
2]{bertrand}. Let $G=A\times L$ be a product of an Abelian variety
by a torus $L$ so that after a finite extension of $k$ %if necessary,
it satisfies Poincar\'e's complete reducibility theorem (as a
variety {\em over $k$}).
 Let  $l$ denote a prime. Let
$G_\linf=\{x\in G(\bar k):\exists n\, l^nx=0\}$ be the
$\linf$-torsion of $G$. For a point $P\in G(k)$, let $G_P$ be
\emph{the smallest algebraic subgroup of $G$ containing $P$}, i.e.
the Zariski closure of subgroup $\Z P$ of $G$, and \emph{let
$G^\circ_{\!P}$ be its connected component through the origin}, and
finally let
\becd\xi_\linf(P):\Gal(\bar k/ k(G_{l^\infty},P)) \longrightarrow  T_{l^\infty}(A\times L)\\
             \sigma \longmapsto \sigma(P_\linf)-P_\linf, \text{ for
             some }P_\linf\in T_\linf(A\times L).
\eecd It can be checked by direct computation the map
$\xi_\linf(P):\Gal(\bar k/ k(G_{l^\infty},P)\ra T_{l^\infty}(A\times
L)$ does not depend on the choice of $P_\linf\in T_\linf(A\times
L)$.

Let $T(G^\circ_{\!P})$ be the sequences of $T(E)$ consisting of
elements of $G^\circ_{\!P}$; then, according to \cite[Theorem
2]{bertrand}, the image of $\xi_\infty(P)=\prod_l \xi_\linf$  has
finite index in $T(G^\circ_{\!P})=\prod_l T_\linf(A\times L)$,
i.e.~the image contains $N T(G^\circ_{\!P})$, for some natural
number $N\in\N$ large enough.

We claim that if we take $G=E^n$, $P=(a_1,\ldots,a_n)\in
E^n(\barq)$, then $N$ above is $N$ required in Lemma. By the result
cited above, it is enough to prove that if $a_1,\ldots,a_n\in
E(\barq)$ are $\EndE$-linearly independent, then
$G^\circ_{\!P}=E^n$.

An algebraic subgroup of $E^n$ is necessarily an Abelian subvariety;
by the Corollary of Poincar\'e's complete reducibility theorem
(Lemma~\ref{lem:poincare:ell} below) if $\dim G^\circ_{\!P} \neq
E^n$, then for some $m\in \N$, the connected component of
$G^\circ_{\!P}/m$ lies in the kernel of some non-trivial morphism
$f:E^n\ra E^m$. The kernel of $f$ satisfies the $\EndE$-linear
relation corresponding to $f:E^n(\bar k)\ra E^n(\bar k)$ as a
morphism of $\EndE$-modules. Since by assumption $P$ satisfies no
$\EndE$-linear relation, this is a contradiction, and
$G^\circ_{\!P}=E^n$.

\subsubsection{Any Abelian subvariety is a connected component of the kernel of a morphism}

The following Lemma has been just used to relate the more geometric
formulation of Bertrand~\cite{bertrand} and the more explicit
statement of Lemma~\ref{lem:kummer:ell} following
Bashmakov~\cite{bashmakov} and
Ribet~\cite{ribet-kummer-theory-of-Abelian-tori} specific to Abelian
varieties.

\begin{lem}[Corollary of Poincar\'e's complete reducibility theorem]\label{lem:poincare:ell}
 Let $B\subset E^n$ be an irreducible Abelian subvariety. Then there
 is a morphism $f:E^n\ra E^m$ such that $B$ is a connected component
 of $\ker f$.

In particular, this implies that an Abelian subvariety $B\neq E^n$
satisfies a non-trivial $\End$-linear relation.
\end{lem}

\bp By Poincar\'e's reducibility theorem, there exist an Abelian
subvariety $B'\subset E^n$ such that $B(\bar k)+B'(\bar k)=E^n(\bar
k)$ and the intersection $B(\bar k)\cap B'(\bar k)$ is finite. This
implies there is an isogeny $f:B\times B' \ra E^n, (x,y)\ra (x+y)$,
and $f(B)=B$. There exists another isogeny $f':E^n\ra B\times B'$
such that $f'\circ f=[m]$ is a multiplication by a natural number
$m\in \N$; let $B_m$ and $B'_m$ be the connected components of
$B/m\subset E^n$ and $B'/m\subset E^n$, respectively, passing
through the origin. Then we know that $m B_m=f f'(B_m)=B$, $m
B'_m=ff'(B'_m)=B'$ but also $f'(B_m)$ and $f'(B'_m)$ do not
intersect within $B\times B'$; to conclude, $B_m,B'_m\subset E^n$
and $B_m\cap B'_m=\emptyset$ and $B_m + B'_m=E^n$. Let $p_2:B\times
B' \ra B'$ be the projection on the second coordinate; since $0$ is
a connected component of $\ker f$, it implies that $B_m$ is a
connected component of kernel of $f\circ p_2 \circ f': E^n\ra E^n$
and $B$ is a connected component of $f\circ p_2 \circ f'\circ [m]:
E^n\ra E^n$.\ep

This concludes the proof of the lemma~\ref{lem:kummer:ell} stating
the Kummer theory for an elliptic curve.\ep

\subsubsection{An inductive argument based on Kummer theory}

 Assume now that we are on inductive step $n-1$, i.e.~we have defined an $\EndE$-linear map
 $h_\nmin:\qV_\nmin\ra W$ and $\sigma_\nmin\in\galqk$, $h_\nmin(v_i)=w_i, 0<i<n$ such that
 $\psi(h_\nmin(v))=\sigma_\nmin\phi(v)$ for every $v\in\qV_\nmin$. Consider a
 compatible system $(\phi(v_1/j))_j,\ldots,(\phi(v_n/j))_j,j\in\N$ of division
 points in $E(\barq)$, and take $N$ as in Kummer theory Lemma~\ref{lem:kummer:ell}. By
 the induction hypothesis we have
 $\sigma_{n-1}\phi(v_1/j)=\psi(v_1/j),\ldots,\sigma_{n-1}\phi(v_{n-1}/j)=\psi(v_{n-1}/j)$ for any $j$.
 Choose $w_n\in W$ such that $\sigma_{n-1}\phi(v_n/N)=\psi(w_n/N)$;
 that is possible by surjectivity of $\psi:W\ra E(\barq)$.
 By Kummer theory lemma, for $N$ large enough, there exists $\sigma'\in\galqk$ such that
 $\sigma'\sigma_{n-1}\phi(v_1/j)=\psi(v_1/j),\ldots,\sigma'\sigma_{n-1}\phi(v_{n-1}/j)=\psi(v_{n-1}/j),
 $ and $\sigma'\sigma_{n-1}\phi(v_n/j)=\psi(w_n/j)$; let
 $\sigma_n=\sigma'\sigma_{n-1}$. By construction we have that
 $\sigma_n|_{\phi(\Q V_{n-1})}=\sigma_{n-1}|_{\phi(\Q V_{n-1})}$ and
 $\sigma_n\phi(v_i/j)=\psi(v_i/j),0<i<n+1$. This implies
 $\sigma_n\phi(v)=\psi(v)$ for arbitrary $v\in\Q V_n$, thereby
 completing the induction step.

After countably many steps we construct a total $\EndE$-linear map
$h=\cup h_n: V\ra W$ and $\sigma:\phi(V)\ra E(\barq)$.  Since
$\phi(V)=E(\barq)$, the Galois map $\sigma$ is defined on the whole
of $E(\barq)$. Since Galois map $\sigma$ is surjective, this implies
$h:V\ra W$ is surjective, too. This completes the proof of
Proposition~\ref{prp:ext}.

The last argument could in fact have been avoided by a little more
careful inductive construction of $h$: instead of always choosing
$w_n$ to match $\phi(v_n)$ we could  have on odd steps pick an
arbitrary $w_n$ and then chosen $v_n$ so that
$\sigma_n(\phi(v_n))=\psi(w_n)$ while on even steps preserving the
old behaviour. It is very easy to force surjectivity of the
constructed map $h:V\ra W$ this way; it is a very common argument in
model theory called ``a back-and-forth argument''. \ep %\ep

\subsection{Concluding Remarks}

\subsubsection{Image of Galois representation}
\begin{rem}\label{rem:chern} Note that for the arguments in \S\ref{sec:galois-image} it is essential
that the image of the Galois action is as large as possible subject
to linear dependencies. However, this is something specific to
elliptic curves and false for higher dimensional Abelian varieties:
one needs to take care of a symplectic form. This observation shows
that the straightforward generalisation to higher dimensional
Abelian varieties is false. See \cite[IV\S7]{gavr-thesis} for a
discussion of a way to do that.
%We discuss this more in Lemma~\ref{lem:kummer}.
\end{rem}

\subsubsection{Kummer theory}

Remark~\ref{rem:chern} says that in higher dimensions,  the base of
the induction breaks down due to additional restrictions on the
image of the Galois action. This does not happen in the the later
steps of the induction process based on Kummer theory.

\begin{rem}[Generalisations of Kummer theory argument]
Since Kummer theory is known in much larger generality, say for a
product of  arbitrary Abelian varieties, complex tori $\Cm$ and
complex lines $\C$ (\cite{bertrand}), it seems
 straightforward
to generalise the Kummer theory argument above to such a product
$A$. Thus, one would prove that if there exists
$h_0:\ker\phi\ra\ker\psi$ and a Galois map $\sigma_0:E(k(T(E))\ra
E(k(T(E)))$ in $\Gal(k(T(E))/k)$ such that $\phi\circ \sigma_0
=h_0\circ \psi$, then there exists $h:V\ra W$ and $\sigma\in\galqk$
making the diagram (\ref{fml:exact:ell:a}) commute. That is, a
morphism between kernel from  $\Lambda\subset V$ to $\Lambda\subset
W$ extends to a morphism on the whole of $V$. Model-theoretically,
this means that the types of the points of universal covering space
lying over algebraic points $A(\barq)$ is \emph{atomic over the
kernel} in the \emph{linear} language.
\end{rem}

%\begin{rem}[A geometric formulations of Kummer theory]
%Note that the formulation of Kummer theory in
%Bertrand~\cite{bertrand} is more geometric than the way we state it
%in Lemma~\ref{lem:kummer:ell}; instead of ``linearly independence of
%coordinates of a point $P\in E^n$'' he speaks of a ``connected
%component of Zariski closure of group $\Z P$ in semi-Abelian variety
%$E^n$'', which is much more invariant and robust notion. The way we
%state Proposition~\ref{asp:atomic} and discussions in
%Chapter~\ref{ch:geom} agree with his formulation well: we also
%prefer to speak only about the properties of \emph{connectivity} and
%\emph{analytic irreducibility}, albeit that of subsets of the
%universal covering space. It is also more close to the functorial
%formulation of Proposition~\ref{prp:functor}.\end{rem}
%
%\begin{rem}[REMORE---------Quantifier elimination and Kummer theory]
%partial\ldotsextends\ldotshomogeneity\ldotsqe!!
%\ldotsluck\ldots\end{rem}

\begin{rem}[Failure of Kummer theory for extensions of Abelian varieties by tori]
According to Ribet \cite{RibetCohom,Ribet}, Kummer theory may fail
for non-trivial extensions of Abelian varieties by $(\Cm)^n$ due to
``existence of an additional morphism''; he gives a motivic
interpretation in \cite{RibetCohom}. It is natural to ask if an
analogous argument could still be carried despite the failure of
Kummer theory. To state a correct conjecture, we may need to use
more general considerations of \cite{gavr-thesis}.
% want to consider
%the functorial formulation of Proposition~\ref{prp:functor} or
%rather that of Proposition~\ref{asp:atomic} in a different language,
%perhaps with respect to \emph{another} variety $A$.
\end{rem}

\section{\label{sec:functor}Fundamental groupoid functor}

In this section we derive Proposition~\ref{prp:functor} from
Proposition~\ref{prp:ext} by carrying out a formal counterpart of a
natural topological construction of a universal covering space.% as a
%space of homotopy classes of paths leaving a basepoint.

\subsection{\label{proof:extTofunc}Uniqueness of extensions implies the universality of the fundamental
groupoid functor}

In this \S\ we show that Proposition~\ref{prp:ext} implies
Proposition~\ref{prp:functor}.  We do so by explicitly constructing
an extension as in Proposition~\ref{prp:ext} from a 2-functor as in
Proposition~\ref{prp:functor}, and then showing that the equivalence
of constructed extensions implies the equivalence of 2-functors.

The construction is a formalisation of the geometric observation
that the universal covering space with a basepoint can be
canonically  identified with the set of homotopy classes of paths
leaving the basepoint of the base.  The identification is via the
unique path-lifting property of the covering map, and depends only
on the choice of a basepoint in the universal covering space. In
case of the universal covering space $\C$ of $E(\C)$, the
correspondence is given by the {\em period map}\nopagebreak
$$\gamma  \longmapsto \int\limits_\gamma dz .$$

\subsubsection{Recovering $V$ as the universal covering space from
the 2-functor $\uOmega$}

Let $\uOmega=(\Pt,\Omega,s,t,\cdot)$ be a 2-functor from $\E$ to
$\Sets$ satisfying the conditions of
Definition~$\ref{functor:cond}$. We want to construct an extension
\becd 0@>>> \Z^2 @>>> V_\uOmega @>\phi>> E(\barq)
@>>> 0\eecd of $\EndE$-modules. %Recall the module $V=V_\uOmega$ is
%an analogue of the universal covering space, and the universal
%covering space can be constructed as the set of all paths leaving a
%particular point, up to homotopy fixing the ends. However, our
%functor $\uOmega$ provides us exactly with the notion of a path up
%to homotopy fixing the ends.
 The topological intuition above suggests that we set
\becd V=V_\uOmega=\bigcup\limits_{y\in \Pt(E)}
\Omega_{0,y}(E)=\{\gamma\in\Omega(E):s(\gamma)=0\} \text{ (disjoint
union)},\\ \phi(\gamma)=t(\gamma)\eecd where $0\in E(k)$ is the zero
point of the elliptic curve $E$.

\subsubsection{Abelian group structure on  $\Omega(E)$}

The functoriality of $\uOmega$ transfers $\EndE$-module structure on
$E(\barq)$ to that on $V$; namely, let us check that the maps
$\Omega(f)$, $f\in\EndE$, and $\Omega(m)$, where $m:E\times E\ra E$
is the morphism of addition on $E$, define $\EndE$-module structure
on $V$, or rather that their restriction to $V$ does.

By assumption the functor $\uOmega$ preserves direct products so
$\Omega(E\times E)=\Omega(E)\times \Omega(E)$, and thus there is a
map
$$ \Omega(m): \Omega(E)\times \Omega(E)\ra \Omega(E).$$
Maps $s,t$ are natural transformations of $\Omega$ to $\Pt$ (as
functors to $\Sets$) and so $s\circ \Omega(m)=\Pt(m)\circ s$,
$t\circ \Omega(m)=\Pt(m)\circ t$ is the map of addition on
end-points. Therefore,
$$\Omega(m)(\Omega_{x,y}(E)\times \Omega_{v,w}(E)) \subset
\Omega_{x+v,y+w}(E)),$$ and in particular
$$\Omega(m)(\Omega_{0,y}(E)\times \Omega_{0,z}(E)) \subset
\Omega_{0,y+z}(E)).$$ Thus $\Omega(m):V\times V\ra V$ gives us a
binary operation. It is straightforward to show that the
preservation of direct product and functoriality implies that
$\Omega(m)$ makes $\Omega(E)$ into an Abelian group. Let us check
this.

By definition, associativity of $m:E\times E\ra E$ means that
$m\circ (m\times\id_E)= m \circ (\id_E \times m):E\times E \times E
\ra E$; by preservation of direct product this implies
$\Omega(m)\circ (\Omega(m)\times\id_E)= \Omega(m) \circ (\id_E
\times \Omega(m))$ and so $\Omega(m)$
 is
associative. Similarly, commutativity of $m$ means
$m\circ(\id_1\times \id_2)=m \circ (\id_2\times \id_1):E\times E\ra
E$; that similarly implies the commutativity of $\Omega(m)$. In the
language of morphisms, the existence of a zero for the additive law
translates to the existence of a morphism $0:\{0\}\ra E$ subject to
the identifies: $m\circ (\id\times 0)=p_2:\{0\}\times E\ra E$ and
$m\circ (0\times \id)=p_1:E\times \{0\}\ra E$ corresponding to a
commutative diagram:
$$\bfig
\morphism(0,0)/<-/<300,300>[E`E\times E;]
\morphism(0,0)/<-/<-300,300>[E`\{0\}\times E;]
\morphism(-300,300)/->/<600,0>[\{0\}\times E`E\times E;]
%[E`E\times E;]
%\morphism(0,600)/<--/<300,-300>[M_3`M_1;]
%\morphism(0,600)/<--/<-300,-300>[M_3`M_2;]
 \efig$$
Apply $\Omega$ to get $\Omega(m)\circ (\Omega(\id)\times
\Omega(0))=\Omega(p_2):\Omega(\{0\})\times \Omega(E)\ra \Omega(E)$
and $\Omega(m)\circ (\Omega(0)\times
\Omega(\id))=\Omega(p_1):\Omega(E)\times \Omega(\{0\})\ra
\Omega(E)$. Preservation of direct product implies $\Omega(m)\circ
(\id\times \Omega(0))=\Omega(p_2):\Omega(\{0\})\times \Omega(E)\ra
\Omega(E)$ and $\Omega(m)\circ (\Omega(0)\times
\id)=\Omega(p_1):\Omega(E)\times \Omega(\{0\})\ra \Omega(E)$. This
implies that $\Omega(0)(\Omega(\{0\}))$ is a zero point in $V$.

Existence of (right) inverse corresponds to the existence of a
morphism $i:E\ra E$ subject to the following commutative diagram:
\becd
E @>>> \{0\} \\
@V(\id_E,i)VV @VV 0 V \\
E\times E @>>m> E \eecd

Again functoriality ensures that $\Omega(i)$ satisfies a similar
diagram, thus proving the existence of inverses.

The above checks that $\Omega(m)$ is an associative commutative
partial operation on $\Omega(E)$ possessing a zero element and
inverses; it is immediate to check $V$ is closed under $\Omega(m)$
and inverse $\Omega(i)$, and so is a group. %, as is indeed $\Omega_{0,0}(E)$.

\subsubsection{Action of ``fundamental group'' $\Omega_{0,0}(E)$ on $V$
via concatenation and by $\Omega(m)$-multiplication}

Take a \emph{loop} $\lambda\in \Omega_{0,0}(E)$ and
$\gamma\in\Omega_{0,y}$; then both concatenation and
$\Omega(m)$-product of $\lambda$ and $\gamma$ are well-defined; let
us show that $\lambda\cdot \gamma=\Omega(m)(\lambda\times\gamma)$:
\[\Omega(m)(\lambda \times \gamma)=\Omega(m)(\lambda\cdot 0
\times 0 \cdot \gamma )= \Omega(m)(\lambda\times 0)\cdot
\Omega(m)(0\times\gamma)=\lambda\cdot\gamma.\]  The latter equality
follows from the inverse element equality of morphisms $m(\id\times
0)=m(0\times \id)=\id$. In the classical example, this observation
corresponds to the following calculation:
$$\int\limits_{\gamma_1\cdot\gamma_2} dz=
\int\limits_{\gamma_1} dz+\int\limits_{\gamma_2} dz=
\int\limits_{\gamma_1\gamma_2} dz.$$ Here ${\gamma_1\cdot\gamma_2}$
denotes the concatenation of the paths and ${\gamma_1\gamma_2}$
denotes the pointwise product of the paths.

\subsubsection{Divisibility of $\EndE$-module structure
and path-lifting property}

Analogously, a morphism $f\in \EndE$, $\Omega(f):\Omega(E)\ra
\Omega(E)$ defines a map $\Omega(f):\Omega(E)\ra \Omega(E) $.
Arguments similar to the ones above allow us to prove that $V$ is an
$\EndE$-module with the operations defined above. Let us now prove
that the path-lifting property implies that $V$ is uniquely
divisible. Indeed, we know that any isogeny $f\in \EndE$ is \'etale
(\cite{MumfordAV}) and we may apply the path-lifting property to get
a bijection
\[\path(f):\bigcup\limits_{0\in\fun(E)}\Omega_{0,y}(E)
\to\bigcup\limits_{z\in\fun(E)}\Omega_{0,z}(E).\] That is, by
 definition of $V$, the map $\Omega(f)$ is a bijection on $V$,
  and $V$ is uniquely divisible as
required.

Finally, to get a short exact sequence as in
Proposition~\ref{prp:ext}, set $\phi(w)=t(w)$. Then naturality of
target map $t$ implies $\phi(w)$ is homomorphism of $\EndE$-modules;
the connectivity of $\Omega(E)$ implies that $\phi$ is surjective.
The kernel $\ker\phi$ of $\phi:V\ra E(\barq)$ is $\Omega_{0,0}(E)$
and is isomorphic to $\Z^2$ by assumption.

\subsubsection{Construction of a natural transformation $h':\uOmega\ra\pitopq\circ \sigma$}
%[{Construction of a natural transformation from ${\underline{\Omega}}$ to $\pitopq$}]

For notational convenience, let
$\uOmega'=(\Pt,\Omega',s',t',\cdot')$ denote the functor $\pitopq$.

Let $W\ra^\psi E(\barq)$ be the $\EndE$-module extension constructed
from $\uOmega'={\pitop}_{|\barq}$ in the same way. Since $V, W$ are
both uniquely divisible extensions of $E(\barq)$ by $\Z^2$, we may
apply Proposition~\ref{prp:ext} to get an $\EndE$-linear map $h:V\ra
W$ and $\sigma\in \galqk$ such that $\sigma\circ\phi=\psi\circ h$.

We want to get a natural transformation $h':\uOmega\ra
\uOmega'\circ\sigma$.  Recall a natural transformation
$h':\uOmega\ra \uOmega'\circ\sigma$ is a family of maps (of sets
with no further structure) $h^\Omega_A:\Omega(A)\ra
\Omega'\circ\sigma(A)$ and $h^\Pt_A: \Pt(A)\ra \Pt'\circ\sigma(A)$
satisfying certain compatibly conditions expressed by commutative
diagremmes
(\ref{fml:nattr:1}),(\ref{fml:nattr:2}),(\ref{fml:nattr:3}).
% and
%the equality $h'(\gamma_1\cdot\gamma_2)=h'(\gamma_1)\cdot
%h'(\gamma_2)$ if the left-hand side product is defined.

Define $h^\Pt_A:A(\barq)\ra A(\barq)$ by $h^\Pt_A(x)=\sigma(x)\in
\Pt'(A)=(\sigma A)(\barq)=A(\barq)$ to be the map induced by Galois
automorphism $\sigma:\barq\ra\barq$.

By assumption of connectivity of $\Omega(E)$,  any $\gamma\in
\Omega(E)$ can be represented as product
$\gamma=\gamma_1^{-1}\cdot\gamma_2$ where $\gamma_1,\gamma_2,
s(\gamma_1)=s(\gamma_2)=0$; define a map $h'_E:\Omega(E)\ra
\Omega'(E)$ by setting $h'_E(\gamma)=h(\gamma_1)\inv\cdot
h(\gamma_2)$. If $\gamma\inv_1\gamma_2=\gamma^{'-1}_1\gamma'_2$ then
$\gamma\inv_1\gamma'_1=\gamma_2\gamma'^{-1}_2\in\Omega_{0,0}(E)$,
and so \becd h(\gamma_1)\inv\cdot h(\gamma_2)=
h((\gamma_1\gamma^{'-1}_1)\gamma'_1)\inv\cdot
h((\gamma_2\gamma^{'-1}_2)\gamma'_2)=\\
h((\gamma\inv_1\gamma'_1)+\gamma'_1)\inv\cdot
h((\gamma_2\gamma^{'-1}_2)+\gamma'_2)= \\
(h(\gamma\inv_1\gamma'_1)+h(\gamma'_1))\inv\cdot
h(\gamma_2\gamma^{'-1}_2)+h(\gamma'_2)=\\
(h(\gamma'_1))\inv\cdot h(\gamma'_2)= h(\gamma^{'-1}_1)\cdot
h(\gamma'_2),
 \eecd
which proves that $h'$ is well-defined. Similar calculations check
that $h'_E$ preserves concatenation $\cdot_E$:
$h'(\gamma\cdot\gamma')=h'(\gamma)\cdot h'(\gamma')$ for arbitrary
$\gamma,\gamma'\in \Omega(E)$.

We define $h'_{E^n}=h^\Omega_{E^n}:\Omega(E^n)\ra\Omega'(E^n)$ from
an arbitrary Cartesian power $\Omega(E^n)=\Omega(E)^n$ by
$h'_{E\times E}(\gamma\times\gamma')=h'_E(\gamma)\times
h'_E(\gamma')$ etc.

%We claim that $h'$ is a natural transformation of 2-functors from
%$\uOmega$ to $\uOmega'\circ\sigma$. We need to check that $h'$ is a
%natural transformation of points functors $\Pt$ to $\Pt\circ\sigma$
%and $\Omega$ to $\pitopq$ compatible with functions $s,t$. This
%follows from the identity $\sigma\circ\phi=\psi\circ h$,
%i.e.~$\sigma\circ t = t \circ h'$.

\comments{Since by definition $h'$ and $\sigma$ coincide on the
functor of points, $h'$ is a natural transformation of points
functors $\Pt_\uOmega$ to $\Pt_\pitopq$. }

\subsubsection{Checking that $h'$ is a natural transformation of $\uOmega$ into
$\uOmega'\circ\sigma$}

To check that $h'$ is a natural transformation of $\uOmega$ to
$\uOmega'\circ \sigma$, we need to check commutativity of the
following diagrams (note that $\sigma$ on the \emph{left}-hand side
is \emph{not} a morphism but a functor!):
\begin{equation}\label{fml:nattr:1}
        \begin{CD} E^n @>>\Pt> E^n(\barq) \\
                  @V \sigma VV   @VV \sigma_{E^n} V \\
                  E^n @> \Pt>> E^n(\barq)\end{CD}
\ \ \ \ \ \ \ \ %
        \begin{CD} A(\barq) @>>\Pt(f) > B(\barq) \\
                   @V \sigma_A VV                  @VV \sigma_B V \\
                   A(\barq) @>\Pt(f\circ\sigma) >> B(\barq)
                   \end{CD}\end{equation}
\begin{equation} \label{fml:nattr:2}       \begin{CD} E^n @>>\Omega> \Omega(E^n) \\
                  @V \sigma VV   @VV h'_{E^n} V \\
                  E^n @> \Omega'>> \Omega'(E^n)\end{CD}
\ \ \ \ \ \ \ \ %
        \begin{CD} \Omega(A) @>>\Omega(f) > \Omega(B) \\
                   @V h'_A VV                  @VV h'_B V \\
                   \Omega'(A) @>\Omega'(f\circ \sigma) >> \Omega'(B) \end{CD}
        \end{equation}
\begin{equation}\label{fml:nattr:3}\begin{CD}
 \Omega(A) @>s_A>>\Pt (A) \\
  @V h'_A VV            @V\sigma_A VV \\
 \Omega'(A) @>s'_A>> \Pt(A)
 \end{CD}\ \ \ \ \ \ \ \ \ %
 \begin{CD}
 \Omega(A) @>>t_A>\Pt (A) \\
  @V h'_A VV            @V\sigma_A VV \\
 \Omega'(A) @>t'_A>> \Pt(A)\end{CD}\end{equation}
\begin{equation}\label{fml:nattr:4}\begin{CD}
 \Omega(A) \times \Omega(A) @D \cdot_A D \mathsf{partial} D \Omega(A) \\
  @V h'_A  \times h'_A VV            @V h'_A VV \\
 \Omega'(A)\times \Omega'(A) @D\cdot'_A D \mathsf{partial} D \Omega'(A)
 \end{CD}\ \ \ \ \ \ \ \ \ %
 \end{equation}

The first pair (\ref{fml:nattr:1}) expresses that
$h^\Pt_A=(\sigma_A)_A$ is a natural transformation of set-valued
functors from $\Pt$ to $\Pt\circ \sigma$; the diagrams are
commutative just by definition of the action of $\sigma\in\galqq$ on
category $\E$. Analogously the second pair (\ref{fml:nattr:2})
expresses that $h'$ is a natural transformation of set-valued
functors $\Omega$ to $\Omega'\circ\sigma$. The first diagram in
(\ref{fml:nattr:2}) says that $h'_{E^n}$ is a map from $\Omega(E^n)$
to $\Omega'(\sigma(E^n))=\Omega'(\sigma(E^n))$; the second one in
(\ref{fml:nattr:2}) expresses the linearity of $h'$ with respect to
the morphisms of $E^n \ra E^m$; that follows from $\EndE$-linearity
of $h:V\ra V$.

The third pair (\ref{fml:nattr:3}) expresses compatibility of the
end-point functions and $h'_A$; for $A=E$, this follows from the
main property $\sigma\circ s =s\circ h$ when restricted to
$W^n\subset\Omega(E^n)$, and that $h'$ preserves concatenation
$\cdot_A$; preservation of direct products allows us to extend this
to arbitrary Cartesian power $E^n$. The diagram (\ref{fml:nattr:4})
expresses the fact that $h'$ preserves concatenation; this is by the
definition of $h'$.

%The fact that $h'$ preserves $s_A,t_A$ is by construction; we also
%remarked before that $h'$ preserves $\cdot_A$.

This concludes the proof that $h'$ is a natural transformation and
that of derivation of  Proposition~\ref{prp:functor} from
Proposition~\ref{prp:ext}. \qed
%
%!!!!!!!!!!!!!!!!!!!!!!!!!!!!!!!!!!!!!!!!!make correct the diagrems!

\section{\label{sec:shimura}Shimura curves}
Arithmetics of Shimura curves is well-studied. In particular, for
Shimura curves, there is a quite explicit description of Galois
action analogous to the results on the Galois action on the Tate
module of an elliptic curve: for curves without complex
multiplication it is a result of Ohta (\cite{ohta}, also
\cite[Theorem 123, p.~90]{clark}); for curves with complex
multiplication this is implied by the explicit description of Galois
group provided by the theory of complex multiplication. We thank
A.Yafaev for pointing and explaining us those results. These results
motivate us to make the following conjecture.

Let $S$ be a connected Shimura curve defined over a number field
$k\subset\barq$, and let  $\SC$ be full subcategory of $\factor
{\text{Var}} {\,\barq}$ consisting of Cartesian powers of finite
\'etale covers of $S$. Assume that $S$ has a $k$-rational point $O$.

Recall we denote $\Pt_{{\pitop}_{|\barq}}(V)=V(\barq)$ and
$\Omega_{{\pitop}_{|\barq}}(V)=\{\gamma\in{\pitop}(V(\C)):
s(\gamma),t(\gamma)\in V(\barq)\}$ is the restriction of ${\pitop}$
to $\barq$-rational points.

\begin{conj}[Universality of fundamental groupoid functor]\label{shimura:functor}
Let $\uOmega=(\Pt,\Omega,s_V,t_V,\cdot_V)$ be a functor from
category $\E$ to $\Groupoids$  such that \bi
\item the functor  of points of  is the functor of $\barq$-rational
points:\\ $\Pt(X)=X(\barq)=\Mor_\SC(O,X)$, $X\in\SC$.
\item
$\uOmega$ preserves direct product:
 $\uOmega(X\times Y)=\uOmega(X)\times\uOmega(Y)$.
\item $\uOmega$  has the unique path-lifting property along \'etale
 morphisms:
 for an \'etale morphism $f:X\ra Y$, a path $\gamma\in\Omega(Y)$
 and a point $x\in \Pt(X)$ such that $\Pt(f)(x)=s(\gamma)$, there exists
 a unique path $\wt\gamma\in \Omega(X)$ such that $\Omega(f)(\wt\gamma)=\gamma$
 and $s(\wt\gamma)=x$.

\ei Assume further that \bi
\item[(4)] $\uOmega(E)$ is a connected groupoid
\item[(5)] there is an isomorphism
$$\Omega_{0,0}(S):=\{\gamma\in \Omega(S): s(\gamma)=t(\gamma)=0\}\cong
\pi_1^\mathrm{top}(S(\C),O).$$\ei
 Then there exists an
automorphism $\sigma\in\galqk$ such that the functors $\uOmega$ and
${\pitop}_{|\barq}\circ\sigma:\E\ra \Groupoids$ are naturally
equivalent:
$$\uOmega\cong {\pitop}_{|\barq}\circ\sigma$$
%\becd \E @>>\uOmega> \Groupoids \\
%   @V\exists\! V\sigma\in\galqk V   @Y\exists\!Y \mathrm{natural\ transformation}Y\\
%            \E @>\pitopq>> \Groupoids \eecd
\end{conj}

It is possible that only a weaker conclusion holds: there exists
{\em finitely many} functors $\F_1,\ldots,\F_n$ satisfying (1)-(5)
such that for any functor $\uOmega$ satisfying (1)-(5) there exists
$\sigma\in\galqk$ such that $\uOmega\circ\sigma$ is naturally
equivalent to one of $\F_1,\ldots,\F_n$.

\bibliographystyle{alpha}
%\bibliography{nieuw,nieuw-ams}
\bibliography{nieuw}

\begin{thebibliography}{Mum70}

\bibitem[Bas72]{bashmakov}
M.~Bashmakov.
\newblock The cohomology of an abelian variety over a number field.
\newblock {\em Russian Mathematical Surveys}, 27:25--70, 1972.

\bibitem[Ber88]{bertrand}
Daniel Bertrand.
\newblock Galois representations and transcendental numbers.
\newblock In {\em New advances in transcendence theory (Durham, 1986)}, pages
  37--55. Cambridge Univ. Press, Cambridge, 1988.

\bibitem[Cla03]{clark}
Pete~L. Clark.
\newblock {\em Local and Global Moduli Spaces of Potentially Quaternionic
  Abelian Surfaces}.
\newblock PhD thesis, Harvard University, April 2003.

\bibitem[Gav05]{gavr-thesis}
Misha Gavrilovich.
\newblock {\em Model Theory of the Universal Covering Spaces of Complex
  Algebraic Varieties}.
\newblock PhD thesis, Oxford University, 2005.
\newblock submitted, available at
  http:/\!\!/misha.uploads.net.ru/misha-thesis.pdf.

\bibitem[Gro]{grothendieck-quillen}
Alexander Grothendieck.
\newblock Pursuing stacks, {\em also known as} {L}ong {L}etter to {Q}uillen.
\newblock http:/\!\!/www.math.jussieu.fr/\!\textasciitilde
  leila/\!mathtexts.php.

\bibitem[JR87]{Ribet}
Olivier Jacquinot and Kenneth~A. Ribet.
\newblock Deficient points on extensions of abelian varieties by {${\bf G}\sb
  m$}.
\newblock {\em J. Number Theory}, 25(2):133--151, 1987.

\bibitem[Lan78]{LangDA}
Serge Lang.
\newblock {\em Elliptic curves: {D}iophantine analysis}, volume 231 of {\em
  Grundlehren der Mathematischen Wissenschaften [Fundamental Principles of
  Mathematical Sciences]}.
\newblock Springer-Verlag, Berlin, 1978.

\bibitem[Lan83]{lang-complex-multiplication}
Serge Lang.
\newblock {\em Complex multiplication}.
\newblock Springer, 1983.

\bibitem[Mum70]{MumfordAV}
David Mumford.
\newblock {\em Abelian varieties}.
\newblock Tata Institute of Fundamental Research Studies in Mathematics, No. 5.
  Published for the Tata Institute of Fundamental Research, Bombay, 1970.

\bibitem[Oht74]{ohta}
M.~Ohta.
\newblock On $l$-adic representations of {G}alois groups obtain from certain
  two-dimensional abelian varieties.
\newblock {\em J. Fac. Soc. Univ. Tokyo Sect. 1A Math.}, 21:299--308, 1974.

\bibitem[Rib79]{ribet-kummer-theory-of-Abelian-tori}
Karl Ribet.
\newblock Kummer theory on extensions of varieties by tori.
\newblock {\em Duke Mathematical Journal}, 46, 1979.

\bibitem[Rib87]{RibetCohom}
Kenneth~A. Ribet.
\newblock Cohomological realization of a family of {$1$}-motives.
\newblock {\em J. Number Theory}, 25(2):152--161, 1987.

\bibitem[Shi71]{Shimura}
G.~Shimura.
\newblock {\em Introduction to the arithmetic theory of automorphic functions}.
\newblock Princton: Prinston University Press, 1971.

\bibitem[VK91]{voevodsky-kapranov}
Vladimir Voevodsky and Mikhail Kapranov.
\newblock $\infty$-groupoids and homotopy types.
\newblock {\em Cah. Top. G\'eom. Diff. Cat.}, 32:29--46, 1991.

\bibitem[Zila]{ZilberCovers}
Boris Zilber.
\newblock Covers of the multiplicative group of an algebraically closed field.
\newblock http:/\!\!/www.maths.ox.ac.uk/\textasciitilde zilber.

\bibitem[Zilb]{ZilberLP}
Boris Zilber.
\newblock Logically perfect structures.
\newblock slides, available at http:/\!\!/www.maths.ox.ac.uk/\textasciitilde
  zilber.

\end{thebibliography}
\end{document}
\subsection{Main ideas}

Consider an elliptic curve $E$ defined over a number field, and
assume for simplicity that $E$ does not have complex multiplication.
Our main observation (Proposition~\ref{prp:ext}) is that {\em the
Galois group $\galqq$ acts transitively on the set of uniquely
divisible abelian group extensions of the group $E(\barq)$ of
algebraic points of an elliptic curve, by $\Z^2$}. Topological
intuition suggests that we think of such an extension $U$ as a
universal covering space of $E(\barq)$, and identify points in $U$
with ``homotopy classes of paths in $E(\barq)$''. This allows to
interpret the observation as saying that the ``Poincar\'e
fundamental groupoid'' of an elliptic curve $E(\barq)$ is uniquely
described by its natural properties
%like unique path-lifting, functoriality etc,
up to a Galois automorphism of $E(\barq)$: the algebraic properties
of the extension translate to the properties of unique path-lifting
and functoriality of the groupoid.  Let us state this  in the
language of category theory. Consider the full subcategory $\E$ of
the category of $\barq$-rational varieties generated by the
Cartesian powers of $E$; the category $\E$ is a subcategory of the
category of complex algebraic varieties. Consider a functor from
$\E$ to groupoids. The observation becomes
Proposition~\ref{prp:functor:0}: {\em Some natural conditions force
such a functor on $\E$
%from the full category of complex algebraic varieties
to be naturally equivalent to the  Poincar\'e fundamental groupoid
functor, perhaps composed with an automorphism of the source
category.} Here we consider the  Poincar\'e fundamental groupoid
functor restricted to the algebraic points. In this form the
question is easily generalised to other varieties, for example for
the class of Shimura curves; however, \cite{gavr-thesis} argues that
the immediate generalisation is not appropriate for many other
situations; see~\cite[Chapter~I]{gavr-thesis} for  generalities and
\cite[\S\S IV.6-7]{gavr-thesis} for more explicit examples.
% This is interesting, as %it may be interpreted as {\em a purely algebraic definition of the
%fundamental groupoid functor}.
\subsection{Arithmetic consequences} The proof of these facts for $E$ depends on the Kummer theory of $E$
and on the image of the Galois action on the torsion points $\Etors$
of $E$; we think that our approach {\em explains the Kummer theory
of $E$}, as well as {\em the results about the image of the Galois
action}. Basically, the proofs show that {\em all the restrictions
on the image of Galois representation on the Tate module of $E$,
come from the geometry of $E$, or the topology of $E$}. Accordingly,
if an analogue of the universality property of the fundamental
groupoid functor holds for Shimura curves, this will provide an
explanation of some results about the arithmetic of Shimura curves.
One may think of Proposition~\ref{prp:functor:0} as a kind of Hasse
principle: for some classes of varieties it is true, and if true, it
provides a {\em formal} way to go from the geometry of a variety, to
(certain questions of) the arithmetics of the variety.

++++++++++++

+++++++++++++++++++++++++++++++++++

Let us now remind the notations related to groupoids. A functor to
the category $\Groupoids$ of groupoids can be represented as a
5-tuple $\uOmega=(\Pt,\Omega,s,t,\cdot)$ consisting of two
set-valued functors $\Pt,\Omega$ where $\Pt(X),\Omega(X)$ is the set
of points and paths(morphisms) of groupoid $\uOmega(X)$,  functions
$s,t$ are \emph{source} and \emph{target} maps giving the end-points
of paths(morphisms), and $\cdot$ is the concatenation of paths
(composition of morphisms).

The category of complex algebraic varieties is naturally embedded in
the category of topological spaces; the embedding allows us to
consider the Poincar\'e fundamental groupoid functor
${\pitop}:\factor {\textrm{Var}}{\,\C}\ra\Groupoids$ as defined on
the category of
 complex algebraic varieties. We consider its restriction %of the functor ${\pitop}$
to the following full subcategory. Let
$$\E\subset\factor{\textrm{Var}}{\,\barq},\ \Ob\, \E=\{E^n:n\geq 0\},\
\Mor_{\!\E}(X,Y)=\Mor_{\text{Var}/\barq}(X,Y)$$ be the full
subcategory of the category of varieties whose objects are the
Cartesian powers of $E$ including the trivial Cartesian power
$E^0=0$. The Galois group $\galqk$ acts on the category
$\text{Var}/\barq$, and its action restricts to the category $\E$:
the action leaves each object of $\E$ invariant and permutes the
morphisms. Note that $\galqk$ acts on the set of functors from $\E$
via its action on
the category $\E$ by automorphisms. %: it acts on the functor of points
%via the formula $\Pt(V)=\sigma(\Pt(V))=\sigma(V(\barq))$ while
%leaving the functor of paths $\Omega $ and functions $s,t$
%invariant.
Recall a groupoid $\Omega(E)$ is \emph{connected} iff for every
$x,y\in\Pt(E)$ there exists $\gamma\in \Omega(E)$ ``going from point
$x$ to point $y$'', i.e.~$x=s(\gamma),y=t(\gamma)$.

The following proposition is derived from Proposition~\ref{prp:ext}
in \ref{proof:extTofunc}.

\begin{prp}[Universality of fundamental groupoid functor]\label{prp:functor:0}\label{prp:functor}
Assume that $E$ defined over a number field $k$ is an elliptic curve
without complex multiplication. Let $\E$ be the category of
Cartesian powers of $E$ as above.

Let $\uOmega=(\Pt,\Omega,s_V,t_V,\cdot_V)$ be a functor from
category $\E$ to $\Groupoids$  such that \bi
\item the functor of points of  is the functor of $\barq$-rational points:
$\Pt(E^n)=E^n(\barq)=\Mor_\E(0,E^n)$, $E^n\in\E$.
\item
$\uOmega$ preserves direct product:
 $\uOmega(X\times Y)=\uOmega(X)\times\uOmega(Y)$.
\item $\uOmega$  has the unique path-lifting property along \'etale
 morphisms:
 for an \'etale morphism $f:X\ra Y$, a path $\gamma\in\Omega(Y)$
 and a point $x\in \Pt(X)$ such that $\Pt(f)(x)=s(\gamma)$, there exists
 a unique path $\wt\gamma\in \Omega(X)$ such that $\Omega(f)(\wt\gamma)=\gamma$
 and $s(\wt\gamma)=x$.

\ei Assume further that \bi
\item[(4)] $\uOmega(E)$ is a connected groupoid
\item[(5)] there is an isomorphism
$$\Omega_{0,0}(E):=\{\gamma\in \Omega(E): s(\gamma)=t(\gamma)=0\}\cong \Z^2
\text{ as }\EndE\text{-modules}.$$\ei
 Then there exists an
automorphism $\sigma\in\galqk$ such that the functors
$\uOmega,{\pitop}_{|\barq}\circ\sigma:\E\ra \Groupoids$ are
naturally equivalent:
$$\uOmega\cong {\pitop}_{|\barq}\circ\sigma$$
%\becd \E @>>\uOmega> \Groupoids \\
%   @V\exists\! V\sigma\in\galqk V   @Y\exists\!Y \mathrm{natural\ transformation}Y\\
%            \E @>\pitopq>> \Groupoids \eecd
\end{prp}

Here $\Pt_{{\pitop}_{|\barq}}(V)=V(\barq)$ and
$\Omega_{{\pitop}_{|\barq}}(V)=\{\gamma\in{\pitop}(V(\C)):
s(\gamma),t(\gamma)\in V(\barq)\}$ is the restriction of ${\pitop}$
to $\barq$-rational points.

!!!!!!!!!!!!!!!!Similarly to Proposition~\ref{prp:ext}, for $E$
possessing complex multiplication, there exists finitely many
$\sigma_i$'s such that any functor on $\E$ as above, is equivalent
to one of $\pitopq\circ\sigma_i$'s.